\documentclass[a4paper,12pt]{article}

\usepackage{xypic,amsmath,amssymb,latexsym,enumerate,color}
\usepackage{theorem,eepic,amsbsy}
 \usepackage{geometry,xspace}
\usepackage[curve]{xy}
\usepackage{float}
\usepackage{epsfig}
\usepackage{epic}
\usepackage{overpic}

\newcommand{\cat}{\mathsf{Cat}}

\def\2cat{\mathrm{2-}\cat}

\def\FTop{\mathsf{FTop}}

\def\crs{\mathsf{Crs}}
\def\CRS{\mathsf{CRS}}

\def\ogpd{$\omega$-$\mathsf{Gpds}$}

\newtheorem{example}{Example}[section]
{\theorembodyfont{\rmfamily}}
{\theorembodyfont{\rmfamily}}
{\theorembodyfont{\rmfamily}}
{\theorembodyfont{\rmfamily}}

\newtheorem{blank}[example]{\hspace{-0.3em}}
{\theorembodyfont{\rmfamily}}

\newtheorem{diag}[example]{Diagram}

\newcommand{\directs}[2]{\def\objectstyle{\scriptstyle} \objectmargin={0pt}
\xy
(0,4)*+{}="a",(0,-2)*+{\rule{0em}{1.5ex}#2}="b",(7,4)*+{\;#1}="c"
\ar@{->} "a";"b" \ar @{->}"a";"c" \endxy }

\newcommand{\xdirects}[2]{\def\objectstyle{\scriptstyle} \objectmargin={0pt}
\xy
(0,0)*+{}="a",(0,-6)*+{\rule{0em}{1.5ex}#2}="b",(7,0)*+{\;#1}="c"
\ar@{->} "a";"b" \ar @{->}"a";"c" \endxy }

\newcommand{\sdirects}[2]{\def\objectstyle{\scriptstyle} \objectmargin={0pt}
\xy
(0,2.2)*+{}="a",(0,-2.5)*+{\rule{0em}{1.5ex}#2}="b",(7,2.2)*+{\;#1}="c"
\ar@{->} "a";"b" \ar @{->}"a";"c" \endxy }

\renewcommand{\ge}{\geqslant}

\def\<{\langle \! \langle}
\def\>{\rangle \! \rangle}

\newcommand{\Ker}{\mbox{Ker} \; }

\def\geq{\geqslant}

\def\B{\mathsf{B}}

\def\Z{\mathbb{Z}}

\def\red{\textcolor{red}}

\def\epsilon{\varepsilon}

\begin{document}

\title{Nonabelian Algebraic Topology}
\author{Ronald Brown\thanks{I am grateful to the IMA for support at this Workshop,
and to the Leverhulme Trust for general support.}} \maketitle
\centerline{\bf UWB Math Preprint 04.15}
\begin{abstract}
This is an extended account of a short presentation with this
title given at the Minneapolis IMA Workshop on `$n$-categories:
foundations and applications', June 7-18, 2004, organised by John
Baez and Peter May.
\end{abstract}

\section*{Introduction}
This talk gave a sketch of a book with the  title {\it Nonabelian
algebraic topology} being written under support of a Leverhulme
Emeritus Fellowship (2002-2004) by the speaker and Rafael Sivera
(Valencia) \cite{rbrsbook}. The aim is to give in one place a full
account of work by R. Brown and P.J. Higgins since the 1970s which
defines and applies crossed complexes and cubical higher homotopy
groupoids. This leads to a distinctive account of that part of
algebraic topology which lies between homology theory and homotopy
theory, and in which the fundamental group and its actions plays
an essential role. The reason for an account at this Workshop on
$n$-categories is that the higher homotopy groupoids defined are
cubical forms of strict multiple categories.

Main applications are to higher dimensional nonabelian methods for
local-to-global problems, as exemplified by van Kampen type
theorems. The potential wider implications of the existence of
such methods is  one of the motivations of this programme.

The aim is  to proceed through the steps
$$\xymatrix{\txt{(geometry)} \ar [r] & \left(\txt{underlying \\ processes}\right) \ar [r] &
\txt{(algebra)} \ar [r] & \txt{(algorithms)} \ar[r] &
\left(\txt{computer \\ implementation} \right) }.$$ The ability to
do specific calculations, if necessary using computers, is seen as
a kind of test of the theory, and one which also leads to seeking
of new results; such calculations seem at this stage of the
subject to require strict algebraic models of homotopy types. We
obtain some nonabelian calculations, and it is this methodology
which we term {\it nonabelian algebraic topology}. It is fortunate
that higher categorical structures, and in particular higher
groupoid structures, do give nonabelian algebraic models of
homotopy types which allow some explicit calculation. They have
also led to new algebraic constructions, such as a nonabelian
tensor product of groups, of Lie algebras, and of other algebraic
structures, with relations to homology of these structures (see
the references in \cite{B-web}).

\section{Background}
Topologists of the early 20th century dreamed of a generalisation
to higher dimensions of the nonabelian fundamental group, for
applications to problems in geometry and analysis for which group
theory had been successful. The dream seemed to be shattered by
the discovery  that \v{C}ech's  apparently natural 1932
generalisation of the fundamental group, the higher homotopy
groups $\pi_n(X,x)$, were abelian in dimensions $\ge 2$. So
`higher dimensional groups' seemed to be just abelian groups, and
the dream seemed to be a mirage.

Despite this, the relative homotopy groups $\pi_n(X,A,x)$ were
found to be in general nonabelian for $n=2$, and as a result
J.H.C. Whitehead in the 1940s introduced the term `crossed module'
for the properties of the boundary map
$$\partial: \pi_2(X,A,x) \to \pi_1(A,x)$$ and the action of
$\pi_1(A,x)$ on $\pi_2(X,A,x)$. In investigating `adding relations
to homotopy groups' he proved the subtle result (call it Theorem
W) that if $X$ is obtained from $A$ by adding 2-cells, then
$\pi_2(X,A,x)$ is the free crossed $\pi_1(A,x)$-module on the
characteristic maps of the 2-cells. The proof used transversality
and knot theory.

A potentially new approach to homotopy theory derived from the
expositions in Brown's 1968 book \cite{B-Book} and Higgins' 1971
book \cite{Higgins4}, which in effect suggested that most of
1-dimensional homotopy theory can be better expressed in terms of
groupoids rather than groups. This led to a search for the uses of
groupoids in higher homotopy theory, and in particular for higher
homotopy groupoids. The basic intuitive concept was generalising
from the usual partial compositions of homotopy classes of paths
to partial compositions of homotopy classes (of some form) of
cubes. One aim was a higher dimensional version of the van Kampen
theorem for the fundamental group. A search for such constructs
proved abortive for some years from 1966.

However in 1974 we observed that Theorem W gave a universal
property for homotopy in dimension 2, which was suggestive. It
also seemed that if the putative higher dimensional groupoid
theory was to be seen as a success it should at least recover
Theorem W. But Theorem W is about relative homotopy groups! It
therefore seemed a good idea for us  to look at the relative
situation in dimension 2, that is, to start from a based pair
$X_*= (X,A,C)$, where $C$ is a {\it set} of base points,  and
define a homotopy double groupoid $\rho(X_*)$ using maps of
squares. The simplest, and most symmetric, possibility seemed to
be to consider $R_2(X_*)$ as the set of maps $f:I^2 \to X$ which
take the edges of $I^2$ to the subspace $A$ and take the vertices
to the set $C$, and then to form $\rho_2(X_*)$ as a quotient by
homotopy {\it rel vertices} of $I^2$ and through elements of
$R_2$. Amazingly, this gave a double groupoid, whose 1-dimensional
part was the fundamental groupoid $\pi_1(A,C)$ on the set $C$. The
proof is not entirely trivial, and uses a filling technique.

Previous to this, work with Chris Spencer in 1971-2 had
investigated the notion of double groupoid, shown a clear relation
to crossed modules, and introduced the notion of {\it connection}
on a double groupoid so as to define a notion of commutative cube
in a double groupoid, generalising the notion of commutative
square in a groupoid. Once the construction of $\rho_2(X_*)$ had
been given, a proof of a 2-dimensional  van Kampen theorem came
fairly quickly. The paper on this was submitted in 1975, and
appeared in 1978 \cite{bh1978} (after a third referee's report,
the reports of the first two referees being withheld). The
equivalence between double groupoids with connection and crossed
modules established by Brown and Spencer \cite{brownspenc:double},
and so relating the double groupoid $\rho_2(X_*)$ to the second
relative homotopy groups, enabled new calculations of some
nonabelian second relative homotopy groups  in terms of the
associated crossed module, and so enabled calculation of some
homotopy 2-types.

Thus there is considerable evidence that the nonabelian
fundamental group is naturally generalised to various forms of
crossed modules or double groupoids, and that these structures
enable new understanding and calculations in dimension 2. The
important point is that in dimension 2 we can easily define
homotopical functors, that is, functors defined in terms of
homotopy classes of maps, in a manner analogous to that of the
fundamental group and groupoid, and then prove, but not so easily,
theorems which enable us to calculate directly and exactly with
these functors without using homological techniques.

As an example, we give  the following problem. Given a morphism
$\iota: P \to Q$ of groups, calculate the homotopy 2-type of the
mapping cone $X$ of the induced map of classifying spaces $B\iota:
BP \to BQ$. This 2-type can be described completely in terms of
the crossed module $\partial: \iota_*(P) \to Q$ `induced' from the
identity crossed module $1: P \to P$ by the morphism $\iota$
\cite{brownwens:ind}.

There are some general results on calculating induced crossed
modules, while for some calculations, we have had to use a
computer. The following table, taken from \cite{brownwens:ind},
contains the results of computer calculations, using the package
XMOD  \cite{XMOD} in the symbolic computation system GAP
\cite{GAP}. The examples are for $Q = S_4$ and various subgroups
$P$ of $Q$. The computer has of course full information on the
morphisms $\partial: \iota_*(P) \to S_4$ in terms of generators of
the groups in the table.    The third column gives $\pi_2 = \Ker
\partial = \pi_2(X)$. The fourth column gives $\pi_1 = \mathrm{Coker}\, \partial = \pi_1(X)$.
The numbers 128.2322 and 96.67 refer to groups
in the GAP4 table of groups.
\begin{center} 
\begin{tabular}{||c|c|c| c ||}
\hline
 $P$ & $\iota_* P$ & $\pi_2 $ &  $\pi_1$      \\
\hline\hline
$ \langle (1,2)\rangle$       & $GL(2,3)$   & $C_2$    & 1     \\
 $S_3$       & $GL(2,3)$   & $C_2$    & 1        \\
$\langle (1,2),(3,4)\rangle$     & $S_4C_2$    & $C_2$  & 1       \\
$D_8$       & $S_4 C_2$  & $C_2$    & 1       \\
 $C_4$       & $96.67$    & $C_4$ & 1         \\
 $C_3$     & $C_3\ SL(2,3)$& $C_6$    &  $C_2$     \\
 $\langle(1,2)(3,4)\rangle $      & $128.2322$   & $C_4C_2^3$  & $S_3$      \\
\hline
\end{tabular}
\end{center}
In fact the pairs of equal groups in the second column give pairs
of isomorphic crossed modules, though the proof of this needs some
calculation.

It is not claimed that these results are in themselves
significant. But they do show that: \\(i) there are feasible
algorithms; (ii) the homotopy groups, even with their structure as
modules over the fundamental group, represent but a pale shadow of
the actual homotopy type; (iii) this homotopy type can sometimes
be represented by convenient nonabelian structures; (iv)  we can
use strict structures for explicit calculations in homotopy
theory; (v) the colimit formulation allows for some complete
calculations, not just up to extension; (vi) there is still the
problem of determining the triviality or not of the first
$k$-invariant in $H^3(\pi_1, \pi_2)$ for the last two examples.

It may be argued that the use of crossed modules has not yet been
extended to the geometric and analytic problems which we have
described above as motivating the need for a higher dimensional
version of the fundamental group and groupoid. This argument shows
that there is work to be done! Indeed, there are few books on
algebraic topology  other than \cite{B-Book} which  mention even
the fundamental groupoid on a set of base points.

Once the definitions and applications of the homotopy double
groupoid of a based pair of spaces had been developed, it was easy
to guess a formulation for a homotopy groupoid and GvKT in all
dimensions, namely replace:\\ the based pair by a filtered space;
\\ the crossed module functor by a crossed complex functor; \\ and
the homotopy double groupoid by a corresponding higher homotopy
groupoid. \\ However the proofs of the natural conjectures proved
not easy, requiring developments in both the algebra of higher
dimensional groupoids and in homotopy theory.

\section{Main results}

Major features of the work over the years with Philip J. Higgins
and others can be summarised in the properties of the following
diagram of categories and functors:

\begin{diag}
 $$\def\labelstyle{\textstyle} \xymatrix@R=3pc{& \txt{\rm filtered spaces}
  \ar @<-0.5ex>[dl] _ -\rho \ar @<0.5ex>[dr] ^ -\Pi& \\
\txt{\rm cubical \\\rm  $\omega$-groupoids \\\rm with connections
} \ar@<-0.5ex>[ur] _-{|\,|\circ U_*}  \ar
@<0.5ex>[rr] ^-\gamma &&\ar@<0.5ex>  [ll]^-\lambda \txt{\rm crossed \\
\rm complexes} \ar @<0.5ex>[ul] ^-{\mathsf{B}}
   }$$
\end{diag}

Let us first say something about the categories  $\FTop$ of
filtered spaces, and \ogpd\ of cubical $\omega$-groupoids with
connections. A {\it filtered space} consists of a (compactly
generated) space $X$ and an increasing sequence of subspaces
$$X_*: X_0 \subseteq X_1 \subseteq X_2 \subseteq
\cdots \subseteq X_n \subseteq \cdots \subseteq   X. $$ With the
obvious morphisms, this yields the category $\FTop$.

One example of a filtered space is a manifold equipped with a
Morse function. This example is studied using crossed complex
techniques in \cite[Chapter VI]{Sharko} (free crossed complexes
are there called homotopy systems).

Another, and more standard,  example of a filtered space is a
$CW$-complex with its skeletal filtration. In particular, the
$n$-cube $I^n$, the product of $n$-copies of the unit interval
$I=[0,1]$, with its standard cell structure, becomes a filtered
space $I^n_*$, and so for a filtered space $X_*$ we can define its
filtered cubical singular complex $R^\Box(X_*)$ which in dimension
$n$ is simply $\FTop(I^n_*,X_*)$.

Thus $R^\Box(X_*)$ is an example of a cubical set: it has  for $ n
\geq 1$ and $ i=1,\ldots,n$,  face maps $\partial^\pm_i $ from
dimension $n$ to dimension $(n-1)$, and in the reverse direction
degeneracy maps $\epsilon _i$, all  satisfying a standard set of
rules. We also need in the same direction as the $\epsilon_i$ an
additional degeneracy structure, called {\it connections},
$\Gamma^\pm_i, \; i=1,\ldots,n$, which are defined by the monoid
structures $\max, \min$ on $I$. There are also  in dimension $n$
and for $i=1,\ldots,n$   partially defined compositions $\circ_i $
given by the usual gluing $a \circ_i b$ of singular $n$-cubes
$a,b$ such that $\partial^+_i a=\partial^-_i b$.

There is an equivalence relation $\equiv$ on $R^\Box(X_*)_n$ given
by {\it  homotopy through filtered maps rel vertices of $I^n$}.
This gives a quotient map \begin{equation} p: R^\Box(X_*) \to
\rho^\Box(X_*)=(R^\Box(X_*)/\equiv).\tag{*}\end{equation} It is
easily seen that $\rho^\Box(X_*)$ inherits the structure of
cubical set with connections. A major theorem is that the
compositions $\circ_i$ are also inherited, so that
$\rho^\Box(X_*)$ becomes a {\it cubical $\omega$-groupoid with
connections}. The proof uses techniques of collapsing on
subcomplexes of the $n$-cube. A development of these techniques to
chains of partial boxes proves the surprising result,  essential
for the theory, that the quotient map $p$ of (*) is a cubical Kan
fibration. This can be seen as a rectification result, by which we
mean a result using deformations to replace equations up to
homotopy by strict equations.

An advantage of the functor $\rho$ and the category \ogpd\ is that
in this context we can prove a Generalised van Kampen Theorem
(GvKT), that {\it the functor $\rho$ preserves certain colimits. }
This yields precise calculations in a way not obtainable with
exact or spectral sequences.

The proof is `elementary' in the sense that it does not involve
homology or simplicial approximation. However the proof is
elaborate. It  depends heavily on the algebraic equivalence of the
categories \ogpd\ and $\crs$ of crossed complexes, by functors
$\gamma, \lambda$ in the Main Diagram. This equivalence enables
the notion of commutative cube in an $\omega$-groupoid, and the
proof that any composition of commutative cubes is commutative. It
is also shown that  $\gamma \rho$ is naturally equivalent to a
functor $\Pi$ defined in a classical manner using the fundamental
groupoid and relative homotopy groups, as detailed below.

In detail we set $\Pi(X_*)_1=\pi_1(X_1,X_0)$, the fundamental
groupoid of $X_1$ on the set of base points $X_0$, and for $n \geq
2$ we let $\Pi(X_*)_n$ be the family of relative homotopy groups
$\pi_n(X_n,X_{n-1},x), \,x \in X_0$. Setting $C_0=X_0$ and
$C_n=\Pi(X_*)_n$ we find that $C$ obtains the structure of {\it
crossed complex},
\begin{diag}
 $$\diagram \cdots \rto & C_n \dto<-.05ex>^{t}\rto^-{\delta_n}
& C_{n-1} \rto \dto<-1.2ex>^{t} & \cdots \rto &
 C_2\rto^-{\delta_2} \dto<-.05ex>^{t} & C_1
\dto<0.0ex>^(0.45){t} \dto<-1ex>_(0.45){s} \\ &
C_0&C_0\rule{0.5em}{0ex}  & & \rule{0.5em}{0ex} C_0 &
\rule{0em}{0ex} C_0  \enddiagram$$
\end{diag}
where the structure and axioms are those universally satisfied by
this example $\Pi(X_*)$.

If $X_*$ is the skeletal filtration of a reduced $CW$-complex $X$,
then $\Pi(X_*)$ should be regarded as a powerful replacement for
the cellular chains of the universal cover of $X$, with operators
from the fundamental group of $X$. However the use of crossed
complexes rather than the cellular chains allows for better
realisation properties. Also the functor $\Pi$ allows for many
base points, so that groupoids are used in an essential way.

The GvKT for $\rho$ immediately transfers to a GvKT for this
classical functor $\Pi$. We again emphasise that the proof of the
GvKT for $\Pi$ is elaborate, and tightly knit, but elementary, in
that it requires no background in homology, or simplicial
approximation. It is a direct generalisation of the proof of the
van Kampen Theorem for the fundamental groupoid. Indeed it was the
intuitions that such a generalisation should exist that gave in
1966 the impetus  to the present theory.

 With this GvKT we deduce in the first instance:
\begin{itemize}
  \item the usual vKT for the fundamental groupoid on a set of
  base points;
  \item  the Brouwer degree theorem  ($S^n$ is $(n-1)$-connected and $\pi _n S^n=\Z$);
  \item   the relative  Hurewicz theorem;
  \item    Whitehead's theorem that  $\pi _2(A \cup
\{e^2_{\lambda} \}_\lambda,A,x)$  is a free crossed
$\pi_1(A,x)$-module;
    \item  a more general excision result on $\pi_n(A \cup B, A,x)$
    as an induced module (crossed module if $n=2$) when  $(A,A\cap
    B)$ is $(n-1)$-connected;
    \item computations of $\pi_n(A \cup B, A \cap B,x)$ when $(A,A
    \cap B),\, (B,A \cap B)$ are $(n-1)$-connected.
  \end{itemize}
  Whitehead's theorem, and the last two results for $n=2$, are nonabelian
  results, and hence not obtainable easily, if at all,  by homological methods.

The assumptions required of the reader are quite small, just some
familiarity with $CW$-complexes. This contrasts with some
expositions of basic homotopy theory, where the proof of say the
relative Hurewicz theorem requires knowledge of singular homology
theory. Of course it is surprising to get this last theorem
without homology, but this is because it is seen as a statement on
the morphism of relative homotopy groups $$\pi_n(X,A,x) \to
\pi_n(X \cup CA,CA,x)\cong \pi_n(X \cup CA,x)$$ and is obtained,
like our proof of Whitehead's theorem , as a special case of an
excision result. The reason for this success is that we use
algebraic structures which model the  geometry  and underlying
processes more closely than those in common use.

The {\it cubical classifying space} $B^\Box C$ construction for a
crossed complex $C$ is given by the geometric realisation of a
cubical nerve $N^\Box C$ defined by
$$(N^\Box C)_n= \crs(\Pi(I^n_*),C).$$ The crossed complex $C$ has
a filtration $C^*$ by successive truncations, and this gives rise
to the filtered space $\B(C)= B^\Box C^*$ for which there is a
natural isomorphism $\Pi\B (C) \cong C$. The fundamental groupoid
 and homotopy groups  of $B^\Box C$ are just the fundamental
 groupoid and `homotopy groups' of $C$. In particular, if $X_*$ is
 the skeletal filtration of a $CW$-complex $X$, and $C= \Pi(X_*)$, then
 $\pi_1(C) \cong \pi_1(X,X_0)$ and for $ n \ge 2, \pi_n(C,x) \cong
 H_n(\widetilde{X}_x)$, the homology of the universal cover of $X$
 based at $x$.

The category $\FTop$ has a monoidal closed structure with tensor
product given by $$(X_* \otimes Y_*)_n= \bigcup _{p+q=n} X_p
\times Y_q.$$ It is also fairly easy to define a monoidal closed
structure on the category \ogpd\,  since this category is founded
 in the structure  on cubes. The equivalence of this category with
that of crossed complexes, $\crs$, then yields a monoidal closed
structure on the latter category, with exponential law
$$ \crs(A \otimes B, C) \cong \crs(A, \CRS(B,C)).$$
The elements of $\CRS(B,C)$ are: in dimension 0 just morphisms $B
\to C$; in dimension 1, are homotopies of morphisms; and in
dimensions $ \geq 2$ are forms of `higher homotopies'. There is a
complicated formula for the tensor product $A \otimes B$, as
generated in dimension $n$ by elements $a \otimes b, \, a \in A_p,
\, b \in B_q, \, p+q =n$, with a set of relations and boundary
formulae  related to the cellular decomposition of a product of
cells $E^p \times E^q$, where $E^p$  has for $p > 1$ just 3 cells
of dimensions $0,p-1,p$ respectively.

Using the filtered 1-cube $I^1_*$, and the values on this of
$\rho$ and $\Pi$, we have unit interval objects in each of the
categories $\FTop$, \ogpd , $\crs$. This, with the tensor product,
gives a cylinder object and so gives a homotopy theory on these
categories.

Another difficult result is that the functor $\Pi$ preserves
certain tensor products, for example of $CW$-filtrations.  This
result is needed to deduce the homotopy addition lemma for a
simplex by induction, and also to prove the homotopy
classification results.

A major result is that there is a natural bijection of homotopy
classes
$$[\Pi(X_*), C] \cong  [X,B^\Box C]$$
for any crossed complex $C$ and $CW$-complex $X$, with its
skeletal filtration. This is a special case of a weak equivalence
$$B^\Box(CRS(\Pi(X_*), C)) \to (B^\Box C)^X.$$ In fact the
published version of this is for  the simplicial classifying space
$B^\Delta C$, and this has given also an equivariant version of
these results.

In this theory the homotopy addition lemma (HAL) for a simplex
falls out by a simple inductive calculation, since the usual
representation of the $n$-simplex $\Delta^n$ as a cone on
$\Delta^{n-1}$ is exactly modelled in the crossed complex case,
where the cone is defined using the tensor product. The explicit
formulae for the tensor product, and so for any cone construction,
give an easy calculation of the HAL formula.

There is an acyclic model theory for crossed complexes similar to
that for chain complexes. So we obtain an Eilenberg-Zilber type
theorem for crossed complexes. This has been developed with
explicit formulae by Tonks.

\section{Conclusion}
The book that is being written is divided into three parts, of
which Part I is now available on the web \cite{rbrsbook}. Part I
deals with the case of dimensions 1 and 2, and the proof and
applications of the GvKT in this dimension. The reasons for this
are that the step from dimension 1 to dimension 2 is a big one,
and the reader needs to grasp the new ideas and techniques. Also a
gentle introduction is needed to calculating with the (nonabelian)
crossed modules.

Part II should be regarded as a kind of handbook on crossed
complexes. The main properties of crossed complexes and of the
functor $\Pi$ are stated. Applications and calculations are
deduced. Thus crossed complexes are presented as a basic tool in
algebraic topology.

The proofs of these basic properties of crossed complexes, such as
the monoidal closed structure, and the GvKT, require however the
use of the category of cubical $\omega$-groupoids with connection.
The theory of these is developed in Part III.

Appendices give some basic theory needed, such as some category
theory (adjoint functors, colimits, etc.).

A question asked at the Workshop was: Why concentrate on filtered
spaces? The only answer I can give is that they provide a workable
basis for a theory of higher homotopy groupoids. A more general
theory, in some ways,  is given by $n$-cubes of spaces, and the
associated cat$^n$-groups, but this does not have such an
intuitive exposition as can be provided for the filtered case.
Currently, there is in the absolute case no theory which works in
 dimensions $>2$, and even the known  dimension 2 case does not yet
provide a range of explicit calculations.

Peter May pointed out that,  in contrast to stable homotopy
theory, this work brings the fundamental group fully into the
algebraic topology setting, and that such input is surely needed
for applications in algebraic geometry. Note  that Grothendieck's
increasingly influential manuscript `Pursuing Stacks' \cite{PS}
was intended as an account of nonabelian homological algebra
(private communication). But that subject should in principle need
`nonabelian algebraic topology' as a precursor.

The book makes no attempt to explain the work on cat$^n$-groups,
which allows for  calculations in homotopy theory which are more
complex, and indeed more nonabelian, than those given in the
planned book. That work thus continues the story of Nonabelian
Algebraic Topology, and  references to many authors and papers are
given in \cite{B-web}.

The references which follow include  some survey articles with
wider references. We conclude with a diagram of historical context
for  the current theory\footnote{I am grateful to Aaron Lauda for
rendering this diagram in xypic. }.

\noindent Professor Emeritus R. Brown,\\
Department of Mathematics, University of Wales, Bangor, Dean St,
Bangor, \\ Gwynedd LL57 1UT,
United Kingdom \\
http://www.bangor.ac.uk/$\sim$mas010  \hspace{1cm} email:
r.brown@bangor.ac.uk

\newpage
\red{\centerline{\bf \large Some Context for Nonabelian Algebraic
Topology }}

\[
 \xy
  (-10,62)*{\xy
        *\txt\sf{Galois \\ Theory}   *+\frm{-,};
        \endxy}="Galois";
  (-35,60)*{ \xy
         *\txt\sf{modulo \\ arithmetic }   *+\frm{-,};
         \endxy}="modulo";
  (-60,60)*{\xy
        *\txt\sf{symmetry }   *+\frm{-,};
        \endxy}="symmetry";
  (20,61)*{ \xy
        *\txt\sf{Gauss' \\ composition of \\ binary quadratic \\forms }   *+\frm{-,};
        \endxy}="Gauss";
  (60,60)*{ \xy
         *\txt\sf{Brandt's \\ composition of \\ quaternary \\ quadratic forms }   *+\frm{-,};
        \endxy}="Brandt";
  (-60,47)*{ \xy
        *\txt\sf{celestial \\ mechanics }   *+\frm{-,};
        \endxy}="celestial";
  (-35,46)*+{ \xy
        *\txt\sf{groups}  *+\frm{-,};
        \endxy}="groups";
  (-8,38)*+{ \xy
        *\txt\sf{van Kampen's \\ Theorem }   *+\frm{-,};
         \endxy}="Van";
  (26,46)*+{\xy
         *\txt\sf{groupoids }    *+\frm{-,};
        \endxy}="groupoids";
  (54,43)*+{ \xy
        *\txt\sf{algebraic \\ topology }    *+\frm{-,};
        \endxy}="algebraic";
  (-67,35)*{\xy
        *\txt\sf{monodromy }    *+\frm{-,};
        \endxy}="monodromy";
  (-39,30)*+{\xy
         *\txt\sf{fundamental \\ group}    *+\frm{-,};
         \endxy}="fundamental";
  (28,30)*{ \xy
        *\txt\sf{invariant\\ theory }    *+\frm{-,};
        \endxy}="invariant";
  (60,25)*{\xy
        *\txt\sf{categories }    *+\frm{-,};
        \endxy}="categories";
  (-63,25)*{ \xy
        *\txt\sf{homology }    *+\frm{-,};
         \endxy}="homology";
(74,34)*{ \xy
        *\txt\sf{homology }    *+\frm{-,};
         \endxy}="homology2";
  (-29,15)*+{\xy
        *\txt\sf{fundamental \\ groupoid }   *+\frm{-,};
        \endxy}="fundgroupoid";
  (3,11)*{ \xy
        *\txt\sf{free \\resolutions}  *+\frm{-,};
        \endxy}="free";
  (40,10)*+{\xy
        *\txt\sf{cohomology \\of  groups } *+\frm{-,};
        \endxy}="cohomology";
  (-62,5)*{\xy
        *\txt\sf{higher \\ homotopy \\ groups \\ (\v{C}ech, 1932) }  *+\frm{-,};
        \endxy}="Cech";
  (-30,-10)*{ \xy
        *\txt\sf{groupoids in \\ differential \\ topology }   *+\frm{-,};
        \endxy}="differential";
  (-5,-5)*+{\xy
        *\txt\sf{identities\\ among\\ relations}   *+\frm{-,};
        \endxy}="identities";
  (24,-10)*+{\xy
        *\txt\sf{double \\ categories}   *+\frm{-,};
        \endxy}="double";
  (60,-4)*{\xy
        *\txt\sf{structured \\ categories \\ (Ehresmann) }   *+\frm{-,};
        \endxy}="Ehresmann";
  (-64,-18)*{\xy
        *\txt\sf{relative \\homotopy\\ groups }  *+\frm{-,};
        \endxy}="relative";
  (-30,-30)*+{\xy
        *\txt\sf{crossed\\ modules }   *+\frm{-,};
        \endxy}="crossed";
  (30,-21)*{\xy
        *\txt\sf{$2$-groupoids }   *+\frm{-,};
        \endxy}="2groupoids";
  (60,-30)*+{\xy
        *\txt\sf{nonabelian \\ cohomology }    *+\frm{-,};
        \endxy}="nonabelian";
  (-60,-40)*{ \xy
        *\txt\sf{crossed\\ complexes}   *+\frm{-,};
         \endxy}="complexes";
  (-6,-28)*+{\xy
        *\txt\sf{double\\ groupoids }   *+\frm{-,};
        \endxy}="doublegroupoid";
  (15,-33)*+{\xy
         *\txt\sf{algebraic \\ $K$-theory}    *+\frm{-,};
         \endxy}="Ktheory";
  (35,-45)*{\xy
        *\txt\sf{ cat$^{1}$-groups \\ (Loday, 1982) }   *+\frm{-,};
         \endxy}="catgroups";
  (-31,-45)*{\xy
        *\txt\sf{cubical \\ $\omega$-groupoids}   *+\frm{-,};
        \endxy}="omega";
  (4,-55)*+{\xy
        *\txt\sf{cat$^{n}$-groups }   *+\frm{-,};
        \endxy}="catn";
  (53,-70)*{\xy
        *\txt\sf{crossed \\ $n$-cubes of \\ groups}    *+\frm{-,};
        \endxy}="ncubes";
  (-65,-56)*{\xy
        *\txt\sf{higher \\ homotopy \\ groupoids}    *+\frm{-,};
        \endxy}="homotopygroupoids";
  (-27,-67)*+{\xy
        *\txt\sf{Generalized \\ van Kampen \\ Theorems}   *+\frm{-,};
         \endxy}="generalized";
  (14,-72)*{\xy
        *\txt\sf{nonabelian \\ tensor \\ products}    *+\frm{-,};
        \endxy}="tensor";
  (-60,-85)*{  \xy
        *\txt\sf{higher order \\symmetry}    *+\frm{-,};
        \endxy}="higherorder";
  (-34,-81)*{\xy
            *\txt\sf{quadratic \\ complexes}   *+\frm{-,};
            \endxy}="quadraticcomplexes";
(-35,-95)*{  \xy
        *\txt\sf{gerbes}    *+\frm{-,};
        \endxy}="gerbes";
  (37,-90)*{\xy
            *\txt\sf{computing \\ homotopy \\ types}   *+\frm{-,};
             \endxy}="computing";
  (-10,-105)*{\xy
             *\txt\sf{multiple \\ groupoids in \\ differential \\ topology}   *+\frm{-,};
             \endxy}="multiple";
(17, -107)*{\xy *\txt\sf{nonabelian\\algebraic\\topology}
            *+\frm{-,};\endxy}="nonabat";
  (60,-105)*{\xy
            *\txt\sf{Pursuing \\Stacks}    *+\frm{-,};
            \endxy}="pursuing";
  {\ar "Galois"; "groups"           };
  {\ar "modulo"; "groups"           };
  {\ar "celestial"; "fundamental"   };
  {\ar "monodromy"; "fundamental"   };
  {\ar "symmetry"; "groups"         };
  {\ar "Gauss"; "Brandt"            };
  {\ar "Gauss"; "groups"            };
  {\ar "groups"; "groupoids"        };
  {\ar "groups"; "fundamental"      };
    {\ar "groups"; "Van"      };
  {\ar "fundamental"; "Van"         };
  {\ar "fundamental"; "fundgroupoid"};
  {\ar "fundamental"; "cohomology"  };
 {\ar "homology2"; "algebraic"  };
 {\ar "Van"; "algebraic"  };
  {\ar "fundamental"; "Cech"        };
 {\ar "doublegroupoid"; "omega"      };
 {\ar  "omega"  ; "homotopygroupoids"      };
  {\ar "Van"; "nonabelian"  \POS?(.75)*{}="Tnonabelian"        };
  {\ar "Van"; "fundgroupoid"        };
  {\ar "fundgroupoid"; "Van"        };
  {\ar "fundgroupoid"; "differential"};
  {\ar "fundgroupoid"; "doublegroupoid"};
  {\ar "fundgroupoid"+(-8,-5); "complexes"};
  {\ar "Brandt"; "groupoids"        };
  {\ar "groupoids"; "categories"    };
  {\ar "groupoids"; "differential"  \POS?(.55)*{}="x" };
  {\ar "x"; "fundgroupoid"  };
 {\ar "algebraic"; "groupoids"    };
  {\ar "algebraic"; "categories"    };
{\ar "homology2"; "categories"    };
  {\ar "algebraic"; "cohomology"    };
  {\ar "free"; "cohomology"         };
  {\ar "free"; "identities"         };
  {\ar "identities"; "crossed"      };
  {\ar "invariant"; "free"          };
  {\ar "homology"; "Cech"           };
  {\ar "Cech"; "relative"           };
  {\ar "categories"; "Ehresmann"    };
  {\ar "Ehresmann"; "double"        };
  {\ar "relative"; "complexes"      };
  {\ar "crossed"; "complexes"       };
  {\ar "complexes"; "omega"         };
  {\ar "complexes"; "homotopygroupoids"};
  {\ar "cohomology"; "Tnonabelian"   };
  {\ar "nonabelian"; "Van"          };
  {\ar "nonabelian"+(6,-6); "pursuing"+(6,6)     };
  {\ar "double"; "doublegroupoid"   };
  {\ar "double"; "2groupoids"       };
  {\ar "doublegroupoid"+(-4,-6); "multiple" \POS?(.14)*{}="x"  \POS?(.46)*{}="y" \POS?(.65)*{}="t"    };
  {\ar "x"; "omega"    };
  {\ar "y"; "generalized"    };
  {\ar "2groupoids"; "doublegroupoid"};
  {\ar "2groupoids"; "Ktheory"      };
  {\ar "Ktheory"; "catgroups"       };
  {\ar "catgroups"+(-12,0); "catn"+(2,3)          };
  {\ar "crossed"; "cohomology"      };
  {\ar "crossed"; "omega"           };
  {\ar "crossed"; "ncubes"   \POS?(.35)*{}="x" \POS?(.75)*{}="z" \POS?(.85)*{}="y"     };
  {\ar "z"; "catgroups"      };
  {\ar "x"; "Ktheory"      };
  {\ar "y"; "nonabelian"      };
  {\ar "crossed"; "doublegroupoid"  };
  {\ar "crossed"; "differential"    };
  {\ar "relative"; "crossed"        };
  {\ar "omega"; "generalized"       };
  {\ar "differential"; "doublegroupoid"      };
  {\ar "catn"; "generalized"        };
  {\ar "homotopygroupoids"; "generalized"};
  {\ar "complexes"; "quadraticcomplexes"+(-7,5)  \POS?(.75)*{}="x" };
  {\ar "x"; "higherorder"   };
  {\ar "quadraticcomplexes"; "computing"};
  {\ar "ncubes"; "computing"};
  {\ar "tensor"; "computing"        };
  {\ar "generalized"; "t"     };
  {\ar "generalized"; "computing" \POS?(.65)*{}="x"    };
  {\ar "x"; "tensor"      };
  {\ar "higherorder"; "pursuing"  \POS?(.25)*{}="x"  };
  {\ar "x"; "multiple"    };
  {\ar "catn"; "ncubes"             };
  {\ar "Ehresmann"; "catgroups"     };
  {\ar "homotopygroupoids"; "catn"     };
  {\ar "catn"; "homotopygroupoids"     };
  {\ar "tensor"+(5,7); "Ktheory"+(4,-6)     };
  {\ar "ncubes"; "tensor"        };
  {\ar "higherorder"; (-60,-120)*{}        };
  {\ar "multiple"; (-10,-120)*{}        };
  {\ar "higherorder"; "gerbes"        };
  {\ar  "gerbes"; "multiple"        };
  {\ar  "gerbes"; (-35,-120)*{}        };
  {\ar "pursuing"; (60,-120)*{}        };
  {\ar "computing"; (37,-120)*{}        };
  {\ar "computing"; "nonabat" };
 {\ar "nonabat"; (17,-120)*{}        };
 \endxy
\]

\end{document}